\documentclass[11pt]{amsart}
\usepackage{amsmath, amsthm, amsfonts, graphicx, overpic}

\newcommand{\R}{\mathbb{R}}
\newcommand{\norm}[1]{\lVert#1\rVert}
\newcommand{\abs}[1]{\lvert#1\rvert}
\newcommand{\Length}[1]{\lvert#1\rvert}
\newcommand{\card}[1]{\lvert#1\rvert}
\newcommand{\epsi}{\varepsilon}
\newcommand{\fhi}{\varphi}
\newcommand{\ipr}[2]{\left\langle #1, #2 \right\rangle}
\DeclareMathOperator{\conv}{conv}
\DeclareMathOperator{\bd}{bd}
\DeclareMathOperator{\interior}{int}

\theoremstyle{plain}
\newtheorem{theorem}{Theorem}
\newtheorem{lemma}{Lemma}
\newtheorem{corollary}{Corollary}
\newtheorem{proposition}{Proposition}
\newtheorem*{lmtheorem}{Lawlor-Morgan Theorem}
\theoremstyle{definition}
\newtheorem{claim}{Claim}

\begin{document}
\title{Three-dimensional antipodal and norm-equilateral sets}

\author{Achill Sch\"urmann}
\address{Department of Mathematics,
        University of Magdeburg,
        39106 Magdeburg, Germany}
\email{\texttt{achill@math.uni-magdeburg.de}}
\author{Konrad J. Swanepoel}
\thanks{This material is based upon work supported by the South African National Research Foundation under Grant number 2053752.}
\address{Department of Mathematical Sciences,
        University of South Africa, PO Box 392,
        Pretoria 0003, South Africa}
\email{\texttt{swanekj@unisa.ac.za}}

\begin{abstract}
We characterize the three-dimensional spaces admitting at least six or at least seven equidistant points.
In particular, we show the existence of $C^\infty$ norms on $\R^3$ admitting six equidistant points, which refutes a conjecture of Lawlor and Morgan (1994, Pacific J. Math \textbf{166}, 55--83), and gives the existence of energy-minimizing cones with six regions for certain uniformly convex norms on $\R^3$.
On the other hand, no differentiable norm on $\R^3$ admits seven equidistant points.
A crucial ingredient in the proof is a classification of all three-dimensional antipodal sets.
We also apply the results to the touching numbers of several three-dimensional convex bodies.
\end{abstract}

\maketitle

\section{Preliminaries}
Let $\conv S$, $\interior S, \bd S$ denote the convex hull, interior and boundary of a subset $S$ of the $n$-dimensional real space $\R^n$.
Define $A+B:=\{a+b:a\in A, b\in B\}$, $\lambda A:=\{\lambda a:a\in A\}$, $A-B:=A+(-1)B$, $x\pm A=A\pm x:=\{x\}\pm A$.
Denote lines and planes by $abc$ and $de$, triangles and segments by $\triangle abc:=\conv\{a,b,c\}$ and $[de]:=\conv\{d,e\}$, and the Euclidean length of $[de]$ by $\Length{de}$.
Denote the Euclidean inner product by $\ipr{\cdot}{\cdot}$.
A \emph{convex body} $C\subset\R^n$ is a compact convex set with nonempty interior.
The \emph{polar} of a convex body $C$ is the convex body $C^\ast := \{x\in\R^n: \ipr{x}{y}\leq 1\text{ for all } y\in C\}$.
Let $\norm{\cdot}$ be a norm on $\R^n$ and denote the resulting normed space, or \emph{Minkowski space}, by $X^n=(\R^n,\norm{\cdot})$.
Denote its unit ball by $B:=\{x:\norm{x}\leq 1\}$.
The \emph{dual norm} $\norm{\cdot}^\ast$ is defined by $\norm{x}^\ast:=\sup\{\ipr{x}{y} : \norm{y}\leq 1\}$.
Denote the dual space by $X^n_\ast=(\R^n,\norm{\cdot}^\ast)$.
Its unit ball is the polar $B^\ast$ of $B$.
See Webster \cite{Webster} for further basic information on convex geometry, and Thompson \cite{MR97f:52001} for the geometry of Minkowski spaces.

\section{Introduction}
An equilateral set $S\subset X^n$ is a set of points satisfying $\norm{x-y}=\lambda$ for all distinct $x,y\in S$, and some fixed $\lambda>0$.
Let $e(X^n)$ be the largest possible size of an equilateral set in $X^n$.
For the Euclidean space $E^n$ with norm $\norm{(x_1,\dots,x_n)}_2=\sqrt{x_1^2+\dots+x_n^2}$ it is a classical fact that $e(E^n)=n+1$.
Petty \cite{MR43:1051} as well as P.~S.~Soltan \cite{MR52:4127} proved that $e(X^n)\leq 2^n$ for all $n$-dimensional normed spaces, and that $e(X^n)=2^n$ if and only if the unit ball is an affine image of an $n$-cube.
Both proved this by showing that equilateral sets are antipodal (see Section~\ref{antipodalsection}), and then using a result of Danzer and Gr\"unbaum \cite{MR25:1488}.
Petty also showed that $e(X^n)\geq 4$ whenever $n\geq3$, and observed that it follows from a result of Gr\"unbaum \cite{MR0159263} on three-dimensional antipodal sets that $e(X^3)\leq 5$ if $X^3$ has a strictly convex norm.
Lawlor and Morgan \cite{MR95i:58051} constructed a smooth, uniformly convex three-dimensional normed space $X^3$ such that $e(X^3)=5$.
Here \emph{smooth} means that the norm is $C^\infty$ on $\R^3\setminus\{o\}$, and \emph{uniformly convex} means that $\norm{\cdot}-\epsi\norm{\cdot}_2$ is still a norm for sufficiently small $\epsi>0$.
They furthermore conjectured \cite[p.~68]{MR95i:58051} that $e(X^3)\leq 5$ for differentiable norms on $\R^3$.
See also Morgan \cite{MR93h:53012}.
Our first result is that this conjecture is false.

\begin{theorem}\label{thm1}
There exists a $C^\infty$ norm on $\R^3$ admitting an equilateral set of six points.
\end{theorem}
Section~\ref{lmsection} provides a simple example, with an equilateral set consisting of a Euclidean equilateral triangle together with a parallel copy rotated by $30^\circ$.
Lawlor and Morgan \cite{MR95i:58051} used equilateral sets to show the existence of certain surface energy-minimizing cones.
In Section~\ref{lmsection} we also describe the cone obtained from the example given in the proof of Theorem~\ref{thm1}.

Proving that $e(X^3)\leq 6$ if the norm is differentiable requires more work.
In particular it involves making a classification of antipodal sets in $\R^3$ (see Section~\ref{antipodalsection}).

\begin{theorem}\label{thm2}
For any differentiable norm on $\R^3$ the size of any equilateral set is at most $6$.
\end{theorem}

Note that by Petty's results \cite{MR43:1051} we have $4\leq e(X^3)\leq 8$, with equality on the right if and only if the unit ball is a parallelepiped.
Along the way in proving Theorem~\ref{thm2} we derive a characterization of the norms admitting at least six or at least seven equilateral points.
The characterization of six equilateral points is in terms of affine regular octahedra and semiregular hexagons.
An \emph{affine regular octahedron with center $o$} is the convex hull of $\{\pm e_1,\pm e_2,\pm e_3\}$, where $e_1,e_2,e_3$ are linearly independent.
Its \emph{one-skeleton} is the union of its $12$ edges.
A \emph{semiregular hexagon $p_1p_2\dots p_6$} is a convex hexagon $\conv\{p_1,p_2,\dots p_6\}$ in some plane of $X^3$ such that all six sides have the same length in the norm, and with $p_1+p_3+p_5=p_2+p_4+p_6$.
In this definition we allow degenerate hexagons where some consecutive sides are collinear.
It is easy to see that a semiregular hexagon of side length $1$ equals $\triangle a_1a_2a_3 - \triangle b_1b_2b_3$ for some two equilateral triangles (in the norm) of side length $1$ in parallel planes.
\begin{theorem}\label{thm3}
Let $X^3$ be a three-dimensional normed space with unit ball $B$, and let $S\subset X^3$ be a set of $6$ points.
Then $S$ is equilateral if and only if
\begin{itemize}
\item either $\bd B$ contains the one-skeleton of an affine regular octahedron $\conv\{\pm e_1,\pm e_2,\pm e_3\}$, and $S$ is homothetic to $\{\pm e_1,\pm e_2,\pm e_3\}$,
\item or $B$ has a two-dimensional face that contains a semiregular hexagon $\triangle a_1a_2a_3-\triangle b_1b_2b_3$ of side length $1$, and $S$ is homothetic to
\[\{a_1,a_2,a_3,b_1,b_2,b_3\}\] where $\triangle a_1a_2a_3$ and $\triangle b_1b_2b_3$ are two equilateral triangles of side length $1$ in parallel planes of $X^3$.
\end{itemize}
In particular, if $S$ is equilateral there always exist two parallel planes each containing three points of $S$.
\end{theorem}
While it may be simple to see if the boundary of the unit ball contains the one-skeleton of an affine regular octahedron (consider for example the rhombic dodecahedron -- Section~\ref{rdex}), it seems to be difficult to determine whether a given $2$-dimensional face contains a semiregular hexagon (Section~\ref{coneex}).
However, by Theorem~\ref{thm3} such faces must have a perimeter of at least $6$, so there cannot be too many of them.

\begin{figure}
\begin{center}
\begin{overpic}[scale=0.27]{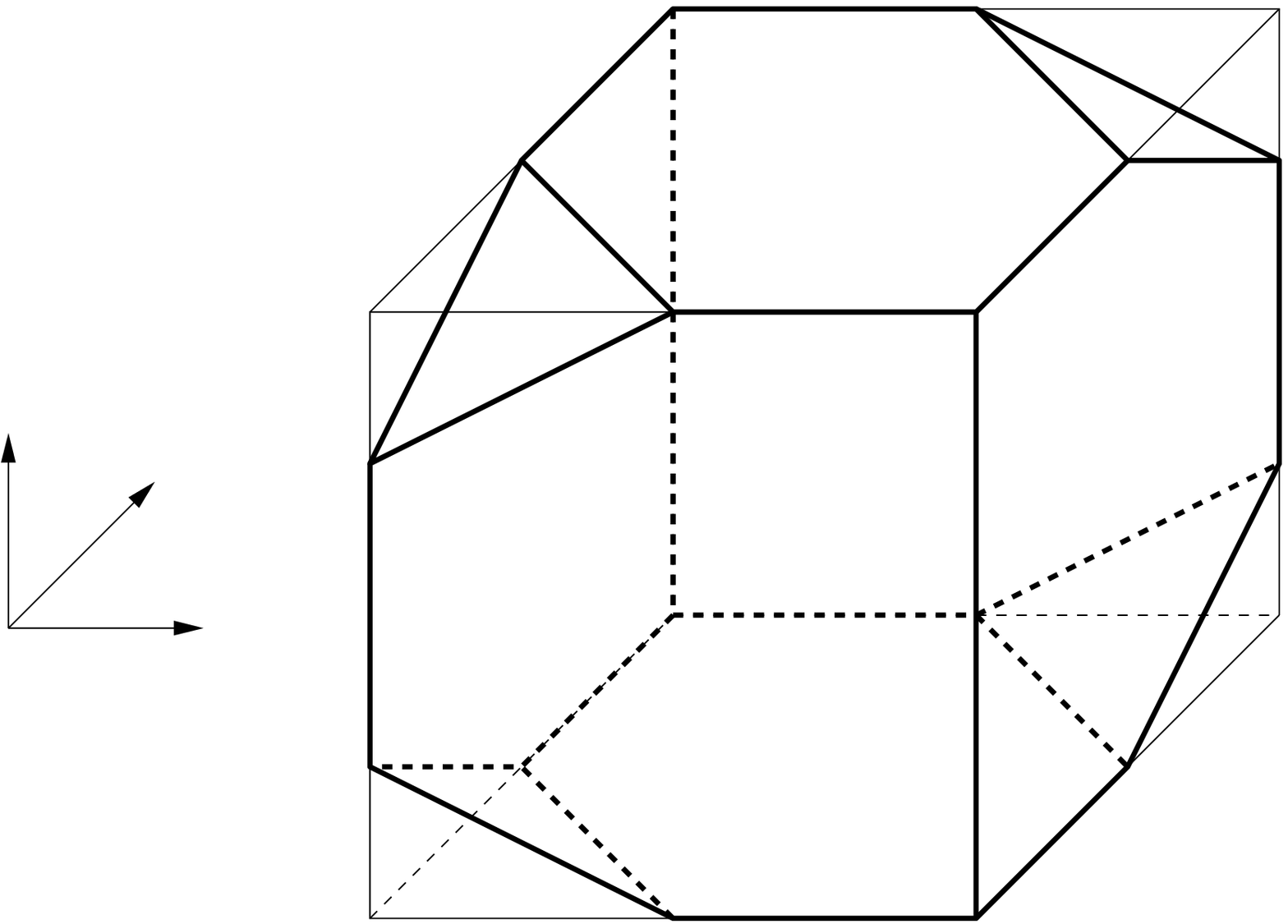}
\put(17,22){$x$}
\put(13,34){$y$}
\put(-1,40){$z$}
\put(20,55){$P_\lambda$}
\end{overpic}
\end{center}
\caption{}\label{fig0}
\end{figure}

The characterization of seven equilateral points is much simpler, as is to be expected.
For $\lambda\in[0,1]$ we define the $3$-polytope $P_\lambda$ to be the polytope with vertex set
\begin{equation*}
\begin{array}{llll}
\pm (-1,1,1), & \pm (1,-1,1), & \pm (-1,0,1), & \pm (1,0,1),\\
\pm (0,1,1), & \pm ( 0,1,-1), & \pm (1,1,-\lambda), & \pm (1,1,1-\lambda).
\end{array}
\end{equation*}
See Figure~\ref{fig0}.
\begin{theorem}\label{thm4}
Let $X^3$ be a three-dimensional normed space with unit ball $B$, and let $S\subset X^3$ be a set of $7$ points.
Then $S$ is equilateral if and only if there exists a linear transformation $\fhi$ and a $\lambda\in[0,1]$ such that $P_\lambda\subseteq \fhi(B)\subseteq [-1,1]^3$, and $\fhi(S)$ is homothetic to
\[ \{(0,0,\lambda), (0,1,0), (0,1,1), (1,0,0), (1,0,1), (1,1,0), (1,1,1)\}.\]
\end{theorem}

Section~\ref{examples} contains applications of Theorems~\ref{thm2}, \ref{thm3} and \ref{thm4}.
Their proofs are given in Section~\ref{lastsection}.

\section{A smooth three-dimensional norm with six equilateral points}\label{lmsection}
\subsection{The construction}\label{constructionsection}
This section does not depend essentially on Theorems~\ref{thm3} or \ref{thm4}.
However, the construction described here is in a sense typical and motivates the development in the remainder of the paper.
We prove Theorem~\ref{thm1} by constructing a $C^\infty$ norm on $\R^3$ admitting an equilateral set of size six.
Note that by Theorem~\ref{thm3}, there will necessarily be two parallel two-dimensional flat pieces on the boundary of the unit ball $B$.
We'll see that $B$ can be chosen such that it has positive curvature at each remaining point of the boundary.

We first construct the equilateral set.
Let $p_k$, $k=0,\dots,11$, be the consecutive vertices of a regular dodecagon $D$ in the $xy$-plane.
To be definite, we may take $p_k=(\cos 2\pi k/12, \sin 2\pi k/12,0)$.
Let $e=(0,0,1)$.
Let $\Delta_1$ be the triangle with vertices $p_3+e$, $p_7+e$, $p_{11}+e$, and $\Delta_2$ the triangle with vertices $p_0$, $p_4$, and $p_8$.
The $\Delta_i$ are congruent equilateral triangles.
We want to construct a smooth norm making $S=\{p_0,p_4,p_8,p_3+e,p_7+e,p_{11}+e\}$ equilateral.
In other words we want to construct a $C^\infty$ unit ball $B$ such that $x-y\in\bd B$ for any two distinct $x,y\in S$.
Let $P=\conv S$.
We first verify that the boundary of $P-P$ contains all $x-y$.
Note that $P-P=\conv(S-S)$, hence $P-P$ is also the convex hull of the union of
\begin{itemize}
\item $\Delta_1-\Delta_2$ in the plane $z=1$,
\item $\Delta_2-\Delta_1$ in the plane $z=-1$,
\item and the regular dodecagon $\sqrt{3}D$ with vertex set
\[ \{\pm(p_0-p_4), \pm(p_0-p_8), \pm(p_4-p_8), \pm(p_3-p_7), \pm(p_3-p_{11}), \pm(p_7-p_{11})\} \]
in the plane $z=0$.
\end{itemize}
Therefore, the hexagons $\pm(\Delta_1-\Delta_2)$ are facets of $P-P$.
It remains to show that the vertices of $\sqrt{3}D$ are all on $\bd(P-P)$.
It is sufficient to show that they are not in the interior of the convex hull $Q$ of the two facets $(\Delta_1-\Delta_2)\cup(\Delta_2-\Delta_1)$.
Note that the intersection of $Q$ with the $xy$-plane is
\[ \frac12(\Delta_1-\Delta_2)+\frac12(\Delta_2-\Delta_1) = \frac12(\Delta_1-\Delta_1)+\frac12(\Delta_2-\Delta_2),\]
which is the dodecagon whose vertices are the midpoints of the edges of $\sqrt{3}D$.
Therefore, the vertices of $\sqrt{3}D$ are on the boundary of $P-P$; even more, they are vertices of $P-P$.
We have shown that each $x-y$, where $x,y$ are distinct points in $S$, is a vertex of $P-P$, except for $\pm(p_7-p_8+e)$, $\pm(p_3-p_4+e)$, $\pm(p_{11}-p_0+e)$, which are in the relative interiors of the facets $\pm(\Delta_1-\Delta_2)$.

It follows that $S$ is equilateral for the norm with unit ball $P-P$.
We now have to smooth $P-P$.
The boundary of any such smoothing $B$ should still contain $\pm(\Delta_1-\Delta_2)$ and the $12$ vertices of $\sqrt{3}D$.
It is well-known that by using convolutions one can construct a $C^\infty$ centrally symmetric convex body $B$ satisfying this requirement --- see e.g.\ \cite[Note~1.3]{Ghomi}.
It follows from the main result of Ghomi \cite{Ghomi} that $B$ can be chosen such that
\begin{itemize}
\item the plane through $\pm(\Delta_1-\Delta_2)$ intersects $B$ in precisely $\pm(\Delta_1-\Delta_2)$,
\item the supporting plane at each vertex $p$ of $\sqrt{3}D$ is perpendicular to the line $op$ (a technical condition needed in Section~\ref{lawlormorgan}), and
\item $\bd B$ has positive curvature everywhere except on $\pm(\Delta_1-\Delta_2)$ and possibly at the $12$ vertices of $\sqrt{3}D$.
\end{itemize}
In fact, by a small modification of the proof in \cite{Ghomi} one can guarantee positive curvature everywhere on $\bd B$ except on $\pm(\Delta_1-\Delta_2)$ (Ghomi, personal communication).
\qed

\subsection{Application to energy-minimizing surfaces}\label{lawlormorgan}
Define the \emph{$\norm{\cdot}$-energy} of a hypersurface $S$ to be $\norm{S}:=\int_S\norm{n(x)}dx$, where $n(x)$ is the Euclidean unit normal at $x\in S$.
Lawlor and Morgan \cite{MR95i:58051} gave a sufficient condition for a certain partition of a convex body by a hypersurface to be energy-minimizing.
We restate a special case of their ``General Norms Theorem~I''.
\begin{lmtheorem}
Let $\norm{\cdot}$ be a norm on $\R^n$, and let $p_1,\dots,p_m\in\R^n$ be equilateral at distance $1$.
Let $\Sigma=\bigcup H_{ij}\subset C$ be a hypersurface which partitions some convex body $C$ into regions $R_1,\dots,R_m$ with $R_i$ and $R_j$ separated by a piece $H_{ij}$ of a hyperplane such that the parallel hyperplane passing through $p_i-p_j$ supports the unit ball $B$ at $p_i-p_j$.

Then for any hypersurface $M=\bigcup M_{ij}$ which also separates the $R_i\cap\bd C$ from each other in $C$, with the regions touching $R_i\cap\bd C$ and $R_j\cap\bd C$ facing each other across $M_{ij}$, we have $\norm{\Sigma}^\ast\leq\norm{M}^\ast$, i.e.\ $\Sigma$ minimizes $\norm{\cdot}^\ast$-energy.
\end{lmtheorem}

Consider the norm $\norm{\cdot}^\ast$ dual to the norm $\norm{\cdot}$ constructed in Section~\ref{constructionsection}.
Since the unit ball $B$ of $\norm{\cdot}$ has two diametrically opposite two-dimensional faces, the dual unit ball $B^\ast$ has two diametrically opposite boundary points $\pm e$ that are not regular --- in fact the set of unit normals of supporting planes at $e$ will be a two-dimensional subset of the Euclidean unit sphere.
Informally, $B^\ast$ is shaped like a spindle.

\begin{figure}
\begin{center}
\includegraphics[scale=0.65]{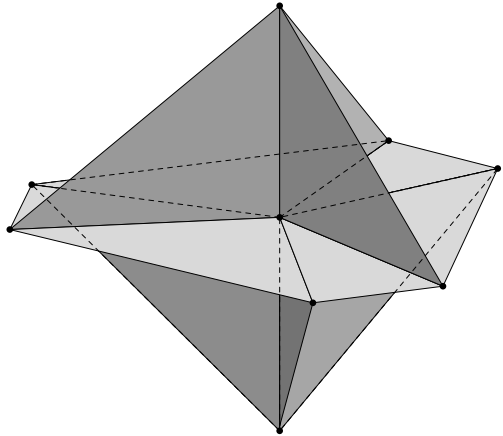}
\end{center}
\caption{Energy-minimizing cone $\Sigma$ with six regions}\label{figms}
\end{figure}

We may now apply the Lawlor-Morgan Theorem as follows.
Consider the equilateral set $S=\{p_0,p_4,p_8,p_3+e,p_7+e,p_{11}+e\}$ of Section~\ref{constructionsection}.
Let $C$ be the convex hull of $\{\pm e, p_2, p_3, p_6, p_7, p_{10}, p_{11}\}$, and let $\Sigma$ be the union of the $12$ triangles
\begin{align*}
&\triangle op_2p_3,\; \triangle op_3p_6,\; \triangle op_6p_7,\; \triangle op_7p_{10},\; \triangle op_{10}p_{11},\; \triangle op_{11}p_2,\\
&\triangle p_3oe,\; \triangle p_7oe,\; \triangle p_{11}oe,\; \triangle p_2o(-e),\; \triangle p_6o(-e),\; \triangle p_{10}o(-e).
\end{align*}
Then $\Sigma$ separates $C$ into six regions (Figure~\ref{figms}).
By the construction of the norm $\norm{\cdot}$ (in particular, the perpendicularity properties), for any $p\in\{p_0,p_4,p_8\}$ and $q\in\{p_3,p_7,p_{11}\}+e$ the supporting plane of $B$ at $p-q$ is parallel to the $xy$-plane, and for any distinct $p,q\in\{p_0,p_4,p_8\}$ or $p,q\in\{p_3,p_7,p_{11}\}+e$, the supporting plane at $p-q$ is perpendicular to $p-q$.
It follows that the hypotheses of the Lawlor-Morgan Theorem are satisfied, giving that $\Sigma$ is $\norm{\cdot}^\ast$-energy-minimizing.

Note that, since $\norm{\cdot}$ is smooth, $\norm{\cdot}^\ast$ is uniformly convex (cf.\ \cite{MR95i:58051}), and since $\bd B$ has positive curvature everywhere except on the two flat pieces, $\bd B^\ast$ is smooth except at $\pm e$.

From Theorem~\ref{thm3} it can be seen that the above example is typical.
For $\norm{\cdot}^\ast$ to be uniformly convex, $\norm{\cdot}$ must be smooth, therefore $B$ must have two-dimensional faces, and then $B^\ast$ must have two nonregular points making $B^\ast$ spindle-shaped.
Because of the two-dimensional faces of $B$ and the structure that the equilateral set necessarily will have, it also follows that the cone $\Sigma$ in the Lawlor-Morgan Theorem must consist of six planar pieces in a plane $\Pi$ parallel to the faces, together with three triangles on one side of $\Pi$ and three triangles on the other side, each with a side on a common line parallel to $oe$.

\section{Applications of Theorems~\ref{thm3} and \ref{thm4}}\label{examples}
\subsection{Regular octahedron}
Bandelt, Chepoi and Laurent \cite{MR1620076} showed that $e(\ell_1^3)=6$, where $\ell_1^3$ is the space with norm $\norm{(\alpha,\beta,\gamma)}_1=\abs{\alpha}+\abs{\beta}+\abs{\gamma}$.
The unit ball is the regular octahedron, and $\{(\pm1,0,0),(0,\pm1,0),(0,0,\pm1)\}$ is clearly equilateral.
To show that $e(\ell_1^3)\leq 6$ using Theorem~\ref{thm4}, it is sufficient to show that no affine regular octahedron contained in $[-1,1]^3$ can contain a $P_\lambda$.
This is easy to see.

\subsection{Rhombic dodecahedron}\label{rdex}
The rhombic dodecahedron $Z$ is the unit ball of the norm $\norm{\cdot}_Z$ with
\[\norm{(\alpha,\beta,\gamma)}_Z:=\max\{|\alpha\pm \beta|,|\alpha\pm \gamma|,|\beta\pm \gamma|\}.\]
The set $\{(\pm1,0,0),(0,\pm1,0),(0,0,\pm1)\}$ is equilateral.
It is again easy to see that no affine rhombic dodecahedron contained in $[-1,1]^3$ can contain a $P_\lambda$.
Therefore, $e(\R^3,\norm{\cdot}_Z)=6$.

\subsection{Spaces and their duals}
As mentioned in the Introduction, for a strictly convex $X^3$ we have $e(X^3)\leq 5$.
The hypothesis of strict convexity cannot be weakened in the following sense.
There exists a unit ball with line segments on its boundary, but no two-dimensional faces, such that $e(X^3)>5$.
Consider for example a ``blown-up octahedron'', where the one-skeleton is fixed (a \emph{wire frame}), but the facets are curved out.
By Theorems~\ref{thm3} and \ref{thm4} we have $e(X^3)=6$ for this norm.
In general we have the following simple consequences of these two theorems.
\begin{corollary}\label{cor}
Let $X^3$ be a three-dimensional normed space.
If the unit ball of $X^3$ does not have a two-dimensional face, then $e(X^3)\leq 6$.
If the unit ball of neither $X^3$ nor its dual has a two-dimensional face, then $e(X^3)\leq 5$.
\end{corollary}

The space $\ell_\infty^3$ has norm $\norm{(\alpha,\beta,\gamma)}_\infty=\max\{\abs{\alpha},\abs{\beta},\abs{\gamma}\}$.
Its unit ball is the cube $[-1,1]^3$, hence $e(\ell_\infty^3)=8$.
Its dual is $\ell_1^3$, for which we know that $e(\ell_1^3)=6$.
Consider now any space $X^3$ with $e(X^3)=7$.
By Theorem~\ref{thm4}, its unit ball $B$ is between some $P_\lambda$ and the cube $[-1,1]^3$.
The polar $B^\ast$ of such a unit ball contains the $1$-skeleton of a regular octahedron on its boundary, and therefore,
$e(X^3_\ast)\geq 6$ by Theorem~\ref{thm3}.
Since $\bd B$ contains an edge of the cube, $\bd B^\ast$ contains two adjacent triangular facets of the octahedron.
It is easily seen that no linear transformation can take $B^\ast$ such that it is between some $P_\lambda$ and $[-1,1]^3$.
By Theorem~\ref{thm4}, $e(X^3_\ast)\leq 6$.
We have shown the following.
\begin{corollary}\label{cor2}
If $e(X^3)\geq 7$, then $e(X^3_\ast)= 6$.
Conversely, if $e(X^3)\leq 5$, then $e(X^3_\ast)\leq 6$.
\end{corollary}

\subsection{Touching numbers}\label{coneex}
Two convex bodies $C,C'\subset\R^n$ \emph{touch} if $C\cap C'\neq\emptyset$ and $\interior C\cap\interior C'=\emptyset$.
For any convex body $C\subset\R^n$ let $C_0:=C-C$ be its \emph{difference body} and let $\norm{\cdot}$ be the norm with unit ball $C_0$, giving a normed space $X^n$.
Let $\{v_1,\dots,v_m\}\subset\R^n$.
The family $\{C+v_i:i=1,\dots,m\}$ is \emph{pairwise touching} if any two translates in the family touch.
It is well known that $\{C+v_i:i=1,\dots,m\}$ is pairwise touching if and only if $\{C_0+2v_i:i=1,\dots,m\}$ is pairwise touching, if and only if $\{v_1,\dots,v_m\}$ is equilateral in $X^n$.
The \emph{touching number} $t(C)$ of $C$ is the largest $m$ such that there exists a pairwise touching family of $m$ translates of $C$.
Then clearly $t(C)=e(X^n)$.
The previous examples show that the touching number of the regular octahedron and the rhombic dodecahedron is $6$.

\begin{figure}
\begin{center}
\includegraphics[scale=0.7]{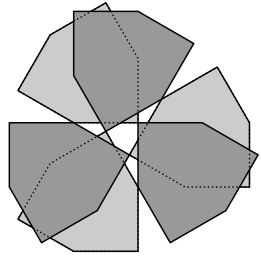}
\end{center}
\caption{A plane supporting six translates of the smooth unit ball $B$}\label{touching}
\end{figure}

The unit ball $B$ of the norm constructed in Section~\ref{constructionsection} has touching number $t(B)=6$.
In particular, there exist six pairwise touching translates of the smooth convex body $B$.
There is a plane, parallel to the $xy$-plane, separating three of the translates from the other three, and with each translate on one side touching each translate on the other side.
This is not easy to visualize and may seem impossible at first.
However, Figure~\ref{touching} shows the intersection of the plane with each translate; there are three translates of the face $\Delta_1-\Delta_2$ touching three translates of the opposite face $\Delta_2-\Delta_1$.
It is easy to see how to modify the construction in Section~\ref{constructionsection} such that $B$ is still smooth but now any pair of the six translates has a two-dimensional intersection.

Consider now any convex disc $D$ in the $xy$-plane of $\R^3$, and let $C$ be the \emph{truncated cone} $\conv(\{e\}\cup D)$, where $e=(0,0,1)$.
For example, if $D$ is a triangle, then $C$ is a tetrahedron and its difference body $C_0$ is the cuboctahedron.
Also, if $D$ is a square, $C$ is a pyramid, and if $D$ is a circular disc, $C$ is a truncated circular cone.
It is easy to see that the touching number of both the tetrahedron and pyramid is at least $5$.
Koolen, Laurent and Schrijver \cite{MR1801196} determined the touching number of the tetrahedron, by showing that $5$ is also an upper bound (see also \cite{MR2034723}).
This is a special case of the following corollary of Theorem~\ref{thm3}.
\begin{corollary}
For any truncated cone $C$ with base a convex disc $D$ we have $t(C)\leq 5$.
\end{corollary}
\begin{proof}
Note that $C-C$ equals the convex hull of $(D-e)\cup(e-D)\cup(D-D)$.
In particular, the extreme points of $C-C$ are contained in the relative boundaries of the three discs $\pm(D-e)$ and $D-D$.
Let $\norm{\cdot}$ be the norm with unit ball $C-C$.
Suppose that $\norm{\cdot}$ has an equilateral set of $6$ points.
Then by Theorem~\ref{thm3} $C-C$ either contains the $1$-skeleton of an affine regular octahedron on its boundary or has a $2$-dimensional face of perimeter at least $6$.

If $\bd(C-C)$ contains the $1$-skeleton of an affine regular octahedron, then the $6$ vertices of the octahedron must be extreme points of $C-C$.
If $D-e$ contains two of these vertices, say $a$ and $b$, then the plane through $\pm a,\pm b$ intersects $C-C$ in the parallelogram with these points as vertices.
In particular, this plane intersects $D-D$ in the segment with endpoints $\pm\frac12(a-b)$.
However, since $[ab]\subset D-e$, it follows that the segment with endpoints $\pm(a-b)$ must be contained in $D-D$, a contradiction.

Therefore, $D-e$, and similarly $e-D$, each contains at most one of the vertices of the octahedron, and it follows that $D-D$ must contain at least $4$ of the vertices.
Therefore, $D-D$ contains exactly $4$ of them, and must be a parallelogram.
Then $D$ is necessarily also a parallelogram, $C$ an affine square pyramid, and $C-C$ the difference body of an affine square pyramid, which is easily seen not to contain the $1$-skeleton of an affine regular octahedron.

In the second case $C-C$ contains a $2$-dimensional face $F$ of perimeter at least $6$.
Suppose $F=D-e$.
It is easy to see that the perimeter of $D-D$ is twice the perimeter of $D-e$ (and more generally, for any two convex bodies $A$ and $B$ in the same plane, the perimeter of $A+B$ equals the sum of the perimeters of $A$ and $B$, in any norm).
The perimeter of $D-D$ is at most $8$, by the theorem of Go\l\c{a}b (see e.g.\ \cite[Theorem~4.3.6]{MR97f:52001}).
Then $D-e$ has perimeter $\leq 4$, a contradiction.

Therefore, $F\neq \pm(D-e)$.
Furthermore, $F$ cannot contain extreme points from both $D-e$ and $-D+e$: if $a+e$ and $-b-e$ are extreme points of $F$, where $a,b\in D$, then their midpoint is $\frac12(a-b)\in\interior(C-C)$, a contradiction.
Thus without loss of generality, the extreme points of $F$ are in $(D-e)\cup(D-D)$.
It follows that $F\cap(D-e)$ and $F\cap(D-D)$ are (possibly degenerate) segments, say $F\cap(D-e)=[ab]-e$ and $F\cap(D-D)=[cd]$, for some $[ab]$ on the relative boundary of $D$ and $[cd]$ on the relative boundary of $D-D$.
(Thus $F$ is either a triangle or a quadrilateral with one pair of opposite edges parallel.)
Without loss, $d-c$ is a positive multiple of $b-a$ if $a\neq b$ and $c\neq d$.
By the definition of $D-D$, $D$ must contain a (possibly degenerate) maximal second edge $[a'b']$ on its relative boundary parallel to $[ab]$ such that 
$a-b'=c$ and $b-a'=d$.
Therefore, $\norm{(a-e)-c}=\norm{b'-e}=1$.
Similarly, $\norm{(b-e)-d}=1$.
Finally, $\norm{(a-e)-(b-e)}=\norm{a-b}\leq\norm{c-d}\leq2$, and it follows that the perimeter of $F$ is at most $6$.
Therefore, it equals $6$, forcing $\norm{a-b}=\norm{c-d}=2$ and $\norm{a'-b'}=0$.
It follows that $\frac{1}{2}(c-d)$ is a unit vector.
Since $\frac{1}{2}(c-d)$ is the midpoint of unit vectors $c$ and $-d$,
all on the relative boundary of $D-D$, the segment $[c,-d]$ is also on the relative boundary.
Therefore, $D-D$ is a parallelogram.
It follows that $D$ is also a parallelogram.
Then $\norm{a'-b'}=\norm{a-b}=2$, a contradiction.

We have shown that neither case in Theorem~\ref{thm3} can occur, and therefore, $t(C)\leq 5$.
\end{proof}

\section{Classifying all antipodal sets in three-space}\label{antipodalsection}
A set $S\subset\R^n$ is \emph{antipodal} if for any two $x,y\in S$ there exist two parallel hyperplanes, one through $x$ and one through $y$, such that $S$ is contained in the closed slab bounded by the two hyperplanes.
See \cite{MS} for a recent survey on antipodal sets.
We recall the following facts.
It is well-known that an antipodal set $S$ is finite, in fact $|S|\leq 2^n$ with equality if and only if $S$ is affinely equivalent to the vertex set of an $n$-cube \cite{MR25:1488}.
It is easily seen that each point of $S$ is a vertex of the polytope $\conv S$.
Two important examples of antipodal sets are equilateral sets in finite-dimensional normed spaces \cite{MR43:1051} (this is how the bound $e(X^n)\leq 2^n$ is deduced) and sets in Euclidean spaces in which no three points span an obtuse angle \cite{MR25:1488}.

In the plane $\R^2$, a set is antipodal if and only if it consists of at most two points, or three noncollinear points, or is the vertex set of a parallelogram.
In $\R^3$, it is clear that any noncoplanar set of four points (the vertex set of a tetrahedron) is antipodal.
By the result of \cite{MR25:1488}, an antipodal set in $\R^3$ has at most $8$ points, with equality if and only if it is the vertex set of a parallelepiped.
In order to characterize three-dimensional antipodal sets it remains to consider sets of size $5$, $6$ and $7$.
Technically the most complicated part is showing that the convex hull of an antipodal set of size $6$ has two parallel facets (Theorem~\ref{sixpoints}).
This has independently been done by Bisztriczky and B\"or\"oczky \cite{BB}.
In fact, they prove this under the weaker requirement that the convex hull is an edge-antipodal polytope.
See also \cite{BBB}.

We constantly refer to the following well-known and easily proved fact.
\begin{lemma}\label{lemma}
A set $S$ is antipodal if and only if for any two distinct points $x,y\in S$, $x-y$ is on the boundary of $\conv(S-S)$.
\end{lemma}
Note that it follows from this lemma that equilateral sets are antipodal, and that an antipodal set $S$ is equilateral in the norm with unit ball $\conv(S-S)$.

\subsection{Five points}
\begin{proposition}\label{fivepoints}
A set of five points in $\R^3$ is antipodal if and only if the points can be labeled as $a, b, c, d, e$ such that $d$ and $e$ are on opposite sides of the plane $abc$, $[de]$ intersects $\triangle abc$ in $p$ such that if we write $p=\lambda a+\mu b+\nu c$ where $\lambda, \mu, \nu\geq 0$, $\lambda+\mu+\nu=1$, then
\begin{equation}\label{star}
\lambda,\mu,\nu\leq\frac{\min\{\Length{dp},\Length{ep}\}}{\Length{de}}. \tag{$\ast$}
\end{equation}
In other words, if we let $\alpha=\min\{\Length{dp},\Length{ep}\}$ and $\beta=\max\{\Length{dp},\Length{ep}\}$, then $p$ must be inside the shaded triangle of Figure~\ref{fig1}.
\begin{figure}
\begin{center}
\begin{overpic}[scale=0.4]{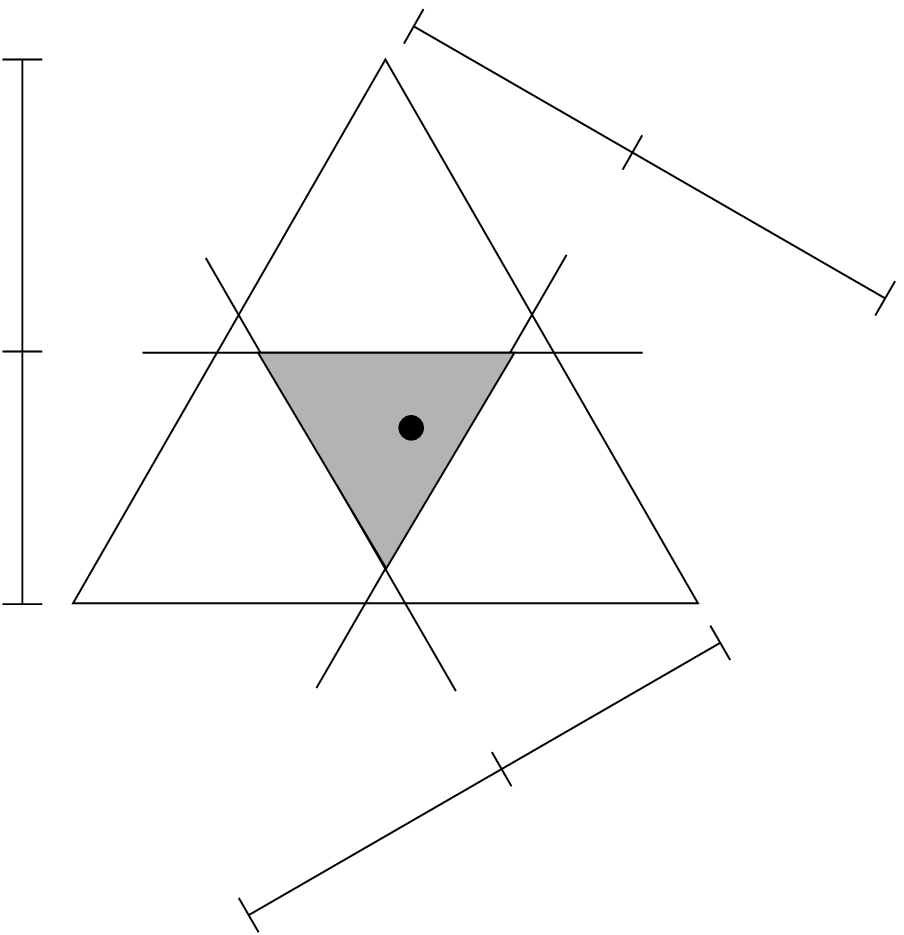}
\put(28,90){$a$}
\put(6,26){$b$}
\put(77,39){$c$}
\put(50,47){$p$}
\put(-7,75){$\beta$}
\put(-7,48){$\alpha$}
\put(37,14){$\beta$}
\put(58,26){$\alpha$}
\put(84,78){$\beta$}
\put(57,92){$\alpha$}
\end{overpic}
\end{center}
\caption{}\label{fig1}
\end{figure}
\end{proposition}

\begin{proof}
Let $S=\{a,b,c,d,e\}$ be antipodal.
Then $\conv S$ has $S$ as vertex set.
It is easily seen, e.g.\ by Radon's theorem, that we may label the points in $S$ such that $[de]$ intersects $\triangle abc$ in a point not in $S$.
Therefore, we may assume without loss of generality in both directions of the proposition that $S=\{a,b,c,d,e\}$ is given so that it is the vertex set of its convex hull, and with $[de]$ intersecting $\triangle abc$ in a point $p=\lambda a+\mu b+\nu c\notin S$ with $\lambda,\mu,\nu\geq 0$ and $\lambda+\mu+\nu=1$.

After applying an appropriate linear transformation we may assume that $\triangle abc$ is equilateral, that $de$ is perpendicular to the plane $abc$, and that $abc$ is parallel to the $xy$-plane (hence $de$ is parallel to the $z$-axis).
Moreover, we may assume that $d$ is in the half space $z<0$ with $\Length{dp}\leq \Length{ep}$.

We show that the nonzero points of $S-S$ are on the boundary of $\conv(S-S)$ if and only if $p$ satisfies \eqref{star}.
Note that $(S-S)\setminus\{o\}$ consists of
\begin{enumerate}
\item the vertices $\{\pm(a-b), \pm(a-c), \pm(b-c)\}$ of a regular
hexagon $H$ in the $xy$-plane, symmetric about $o$,
\item the vertices $\{a, b, c\}-d$
of a triangle $\Delta$ in the half space $z>0$, and the vertices 
of its negative $-\Delta$ in $z<0$,
\item the vertices $e-\{a,b,c\}$
of a triangle $\nabla$ in the half space $z>0$, and the vertices 
of its negative $-\nabla$ in $z<0$,
\item the point $e-d$ in $z>0$, and $d-e$ in $z<0$.
\end{enumerate}
Since $p\in\conv\{a,b,c\}$, it follows that if we orthogonally project the part of $S-S$ in the half space $z\geq 0$ onto the $xy$-plane, we obtain the situation in Figure~\ref{fig2}.
\begin{figure}
\begin{center}
\begin{overpic}[scale=0.4]{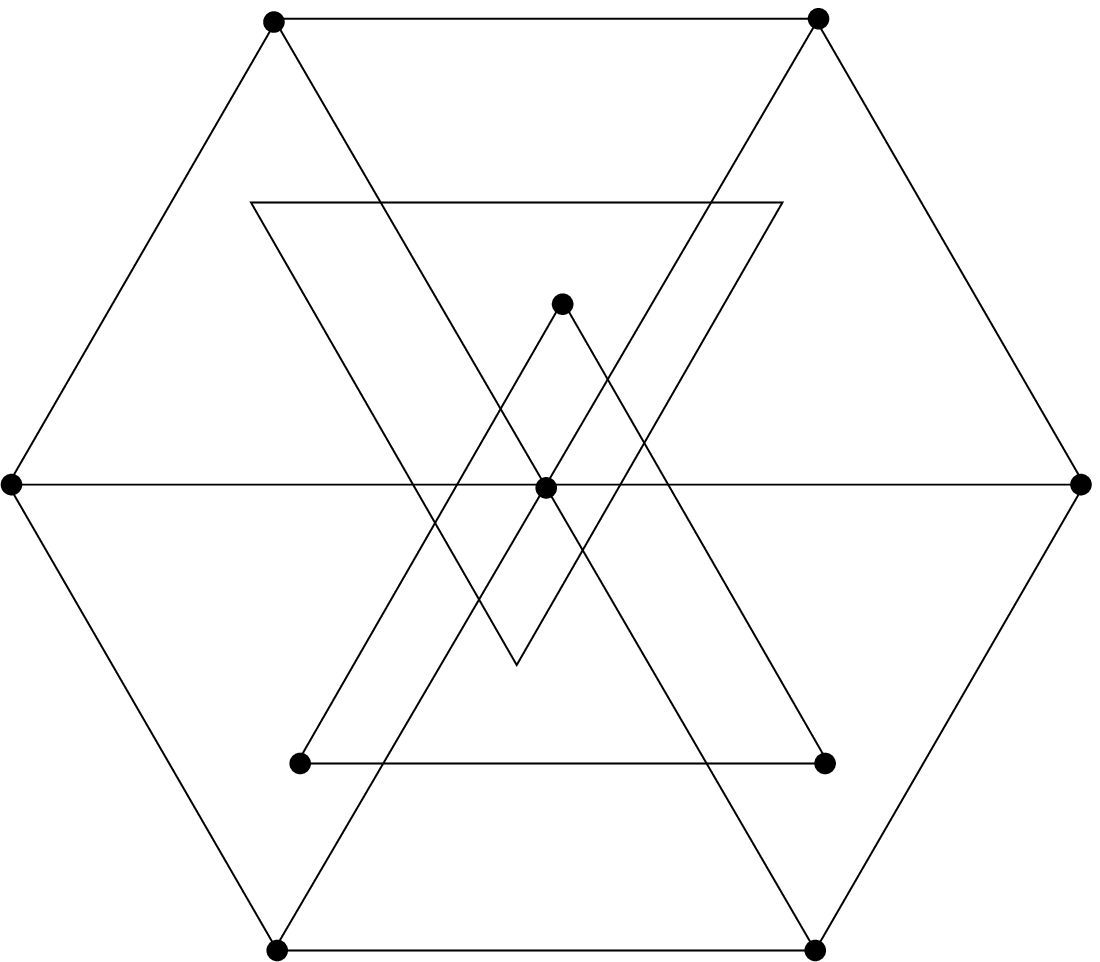}
\put(10,0){$\scriptstyle b-a$}
\put(79,0){$\scriptstyle c-a$}
\put(-14,43){$\scriptstyle b-c$}
\put(53,45){$\scriptstyle o$}
\put(103,43){$\scriptstyle c-b$}
\put(8,85){$\scriptstyle a-c$}
\put(80,85){$\scriptstyle a-b$}
\put(15,20){$\scriptstyle b-p$}
\put(74.5,22){$\scriptstyle c-p$}
\put(45,63){$\scriptstyle a-p$}
\end{overpic}
\end{center}
\caption{}\label{fig2}
\end{figure}
Since a similar picture holds for the part of $S-S$ in $z\leq 0$, it follows that $(S-S)\setminus\{o\}$ is on the boundary of $\conv(S-S)$ if and only if $e-d$ and the vertices of $\Delta$, $\nabla$ and $H$ are on the boundary of $\conv(\{e-d\}\cup\Delta\cup\nabla\cup H)$, i.e., we only have to consider the upper half plane $z\geq 0$.
Clearly $e-d$ and $H$ will be on the boundary.
It remains to show that the vertices of $\Delta$ and $\nabla$ are on the boundary if and only if $p$ satisfies \eqref{star}.
We first show
\begin{claim}\label{claim1}
The vertices of $\Delta$ are not in the interior of the truncated cone $\Gamma=\conv(\{e-d\}\cup H)$, if and only if $p$ satisfies \eqref{star}.
\end{claim}
With $\alpha=\Length{dp}$ and $\beta=\Length{ep}$ we know that $e-d$ is in the plane $z=\alpha+\beta$, and $\Delta$ is in the plane $z=\alpha$.
By projecting the slice $z=\alpha$ of $\Gamma$ onto the $xy$-plane, we see that no vertex of $\Delta$ is in $\interior\Gamma$ if and only if no vertex of the projection of $\Delta$ is in the interior of the hexagon $\frac{\beta}{\alpha+\beta}H$.
See Figure~\ref{fig3}.
\begin{figure}
\begin{center}
\begin{overpic}[scale=0.55, viewport=0 65 313 273, clip=true]{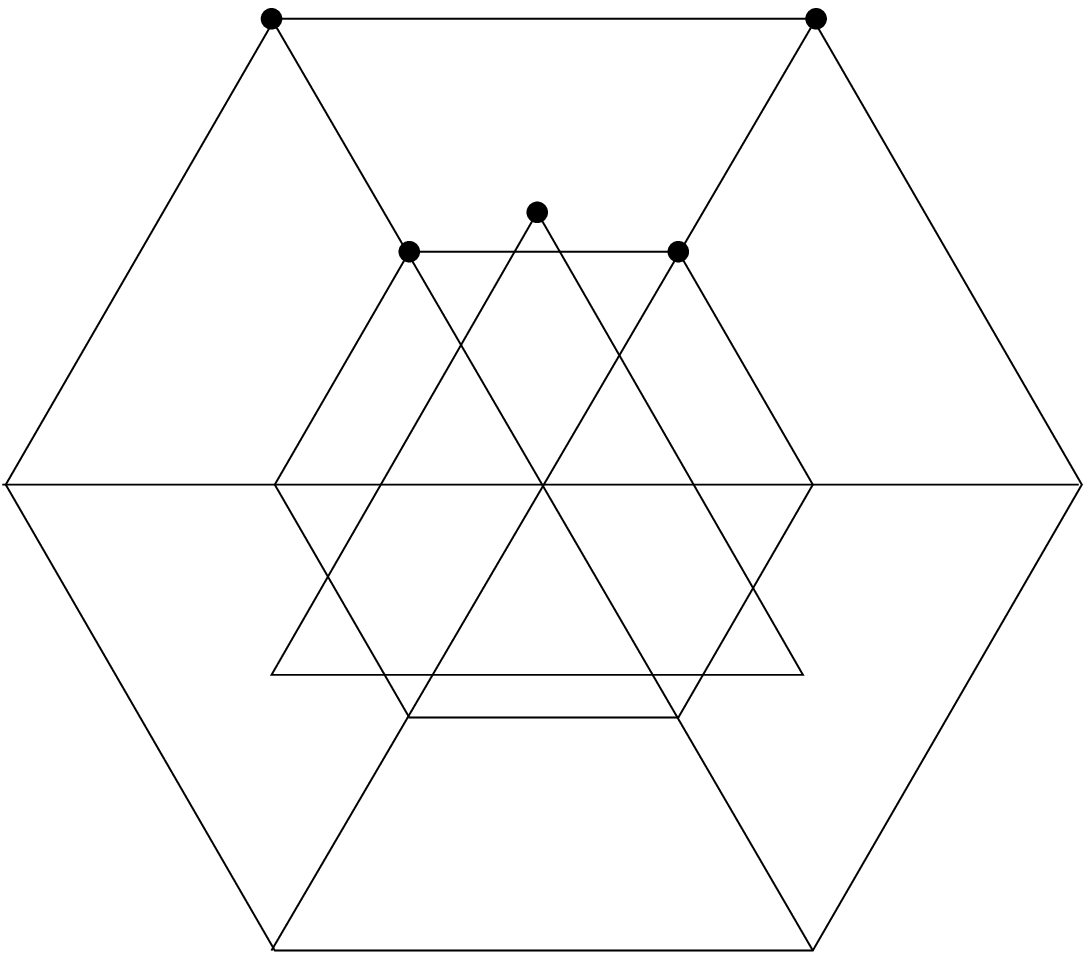}
\put(54,24){$\scriptstyle o$}
\put(13,65){$\scriptstyle a-c$}
\put(78,65){$\scriptstyle a-b$}
\put(13,43){$\scriptstyle \frac{\beta}{\alpha+\beta}(a-c)$}
\put(45,50){$\scriptstyle a-p$}
\put(64,43){$\scriptstyle \frac{\beta}{\alpha+\beta}(a-b)$}
\end{overpic}
\end{center}
\caption{}\label{fig3}
\end{figure}
The projection $a-p$ of $a-d$ is in the triangle $a-\triangle abc$.
Then $a-p\notin\interior\frac{\beta}{\alpha+\beta}H$ if and only if $a-p$ and $o$ are not in the same open half plane of the $xy$-plane bounded by the line through $\frac{\beta}{\alpha+\beta}(a-b)$ and $\frac{\beta}{\alpha+\beta}(a-c)$.
This is easily seen to be equivalent to $\lambda\leq\frac{\alpha}{\alpha+\beta}$.
Similar considerations for $b-d$ and $c-d$ establish Claim~\ref{claim1}.

Since $\nabla$ has a larger $z$-coordinate than $\Delta$ (from $\Length{dp}\leq\Length{ep}$), and the projections of $\nabla$ and $\Delta$ are reflections in $o$, it follows that if $\Delta$ is outside $\interior\Gamma$, then $\nabla$ is also outside $\interior\Gamma$.
It then remains to show
\begin{claim}\label{claim2}
The vertices of $\Delta$ are not in the interior of the half cuboctahedron $\Sigma=\conv(\nabla\cup H)$ if and only if $P$ satisfies \eqref{star}.
See Figure~\ref{fig4}.
\begin{figure}
\begin{center}
\begin{overpic}[scale=0.4]{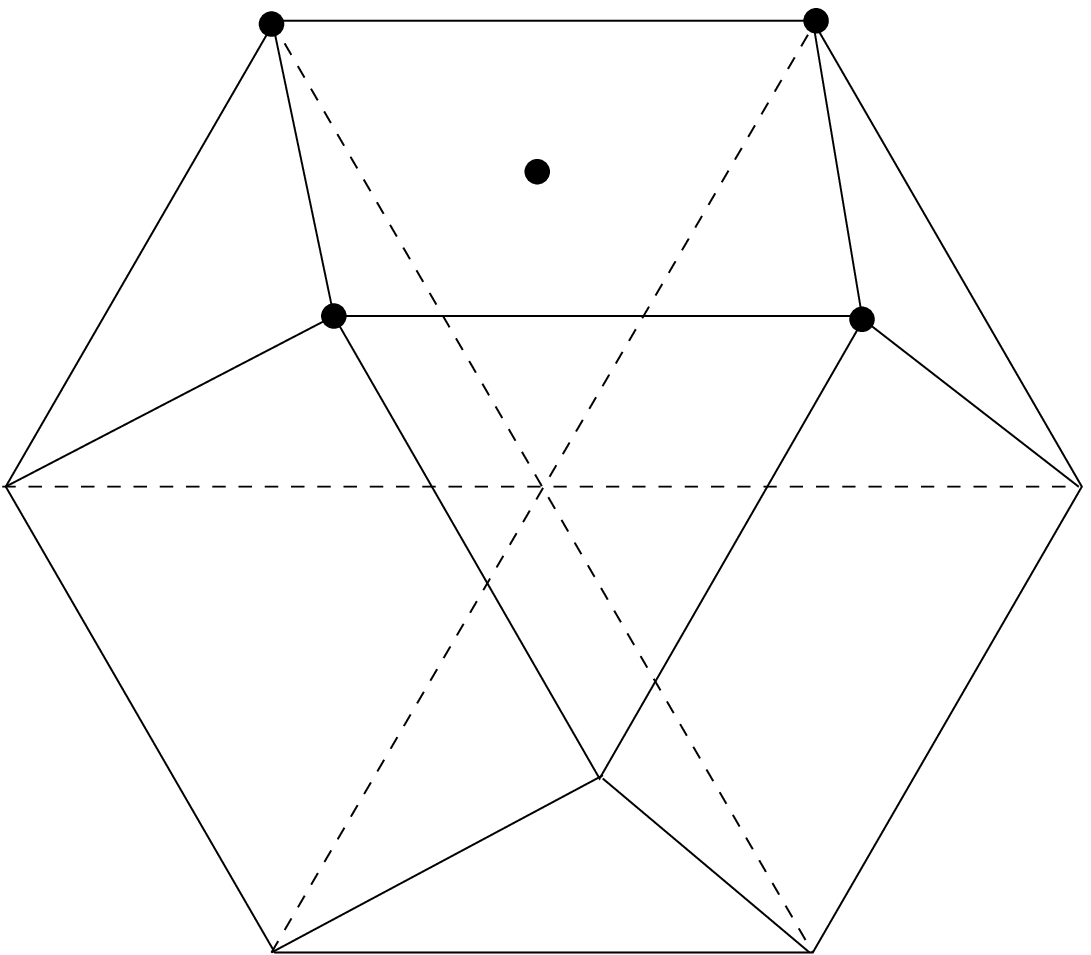}
\put(8,86){$\scriptstyle a-c$}
\put(80,86){$\scriptstyle a-b$}
\put(15,59){$\scriptstyle e-c$}
\put(43,76){$\scriptstyle a-d$}
\put(65,62){$\scriptstyle e-b$}
\end{overpic}
\end{center}
\caption{}\label{fig4}
\end{figure}
\end{claim}
Note that $a-d$ is outside $\interior \Sigma$ if and only if $a-d$ and $o$ are not in the same half space bounded by the plane through the parallelogram with vertex set $\{a-c, a-b, e-c, e-b\}$.
Also, $a-d$ is in the plane $z=\alpha$, which intersects the parallelogram in the line through $\frac\alpha\beta(e-c)+(1-\frac\alpha\beta)(a-c)$ and $\frac\alpha\beta(e-b)+(1-\frac\alpha\beta)(a-b)$.
Projecting onto the $xy$-plane, we find that $a-d\notin\interior\Sigma$ if and only if $a-p$ and $o$ are not in the same open half plane bounded by the line through $\frac\alpha\beta(p-c)+(1-\frac\alpha\beta)(a-c)$ and $\frac\alpha\beta(p-b)+(1-\frac\alpha\beta)(a-b)$, which is easily seen to be equivalent to $\lambda\leq\frac{\alpha}{\alpha+\beta}$.
Similar considerations for $b-d$ and $c-d$ then give Claim~\ref{claim2}.
\end{proof}

\subsection{Six points}
In the sequel we only need the following two consequences of Proposition~\ref{fivepoints}.

\begin{lemma}\label{fivepointscor}
Let $S=\{a, b, c, d, e\}\subset\R^3$ be an antipodal set such that $[de]$ intersects $\interior\conv S$.
Then the following planes support $\conv S$:
\begin{enumerate}
\item the plane through $e$ that contains lines parallel to $ab$ and $cd$,\label{first}
\item the plane through $a$ parallel to $bcd$.\label{second}
\end{enumerate}
\end{lemma}

\begin{proof}
\eqref{first}.
Consider the plane through $ab$ that contains a line parallel to $cd$.
Let $e'$ be the intersection of this plane with $de$.
Note that it is sufficient to prove that $e'\in[de]$.
Let $de$ intersect $\triangle abc$ in $p$.
Let the line through $p$ parallel to $ab$ intersect $ac$ in $q$, and let $cp$ intersect $ab$ in $r$.
Then similar triangles give $\Length{e'p}/\Length{pd}=\Length{rp}/\Length{pc}=\Length{aq}/\Length{qc}$.
By \eqref{star} we must have
\[ \Length{aq}/\Length{qc}\leq\min\{\Length{ep},\Length{pd}\}/\max\{\Length{ep},\Length{pd}\}\leq \Length{ep}/\Length{pd}.\]
It follows that $\Length{ep}\geq \Length{e'p}$, as required.

\eqref{second}.
By the first part of this lemma, the plane through $d$ containing lines parallel to $ae$ and $bc$ supports $\conv S$.
It follows that the plane $bcd$ separates $\conv S$ from the ray emanating from $d$ in the direction $e-a$.
Translating $bcd$ so that it passes through $a$ we obtain that the ray from $a$ through $e$ and the points $b,c,d$ are on the same side of the translated plane, i.e., it supports $\conv S$ at $a$.
\end{proof}

The next proposition describes a construction of antipodal sets of six points which generalizes the construction in Section~\ref{constructionsection}.

\begin{proposition}\label{construction}
Let $\Pi_a$ and $\Pi_b$ be two parallel planes in $\R^3$.
Let $a_1,a_2,a_3\in\Pi_a$ and $b_1,b_2,b_3\in\Pi_b$.
Then the following are equivalent:
\begin{enumerate}
\item \label{one}
The set $S=\{a_1,a_2,a_3,b_1,b_2,b_3\}$ is antipodal.
\item \label{two}
None of the $12$ (not necessarily distinct) points $a_i-a_j$, $b_i-b_j$, $i\neq j$, is in the relative interior of the convex hull of the remaining $11$.
\end{enumerate}
\end{proposition}
\begin{proof}
By Lemma~\ref{lemma} we have to show that \eqref{two} is necessary and sufficient for the nonzero points in $S-S$ to be on the boundary of $D:=\conv(S-S)$.
The points $a_i-b_j\in\Pi_a-\Pi_b$ and $b_i-a_j\in\Pi_b-\Pi_a$ are all clearly on $\bd D$, in the facets $\pm F:=D\cap\pm(\Pi_a-\Pi_b)$.
Therefore, we only have to consider the $12$ points $a_i-a_j$, $b_i-b_j$, $i\neq j$.
Condition~\eqref{two} is clearly necessary for them to be on the boundary.
To see that \eqref{two} is also sufficient, we only have to show that the section $\Sigma$ of $\conv(F\cup-F)$ by the plane through the origin parallel to $\Pi_a$ and $\Pi_b$ is contained in the polygon $P$ with vertex set $a_i-a_j$, $b_i-b_j$, $i\neq j$.
This follows upon noting that
\[ \Sigma = \conv\{\frac12(a_i-b_j)+\frac12(b_k-a_\ell) : 1\leq i,j,k,\ell\leq 3\},\]
and
\[ \frac12(a_i-b_j)+\frac12(b_k-a_\ell)=\frac12(a_i-a_\ell)+\frac12(b_k-b_j)\in P.\qedhere\]
\end{proof}

In the next theorem we show that any $6$-point antipodal set in $\R^3$ is as described in the above proposition.
We also describe all the combinatorial types of their convex hulls.

\begin{theorem}\label{sixpoints}
Let $S$ be an antipodal set of $6$ points in $\R^3$.
Then there exist two parallel planes $\Pi_1$ and $\Pi_2$ such that $\card{S\cap\Pi_i}=3$, $i=1,2$ (thus $S$ is as in Proposition~\ref{construction}).
Furthermore, $\conv S$ is of one of the following two types:
\begin{enumerate}
\item combinatorially equivalent to an octahedron, with some two opposite facets parallel,

\item a ``skew'' triangular prism with one facet a parallelogram with vertices $\{a,b,c,d\}$, and an edge $[ef]$ which is a translate of some segment $[e'f']$ where $e'\in[ad]$ and $f'\in[bc]$ (hence $ade$ and $bcf$ are parallel planes).
There are two combinatorial types, depending on whether $ef$ is parallel to $ab$ or not.
See Figure~\ref{fig5}.

\begin{figure}
\begin{center}
\begin{overpic}[scale=0.4]{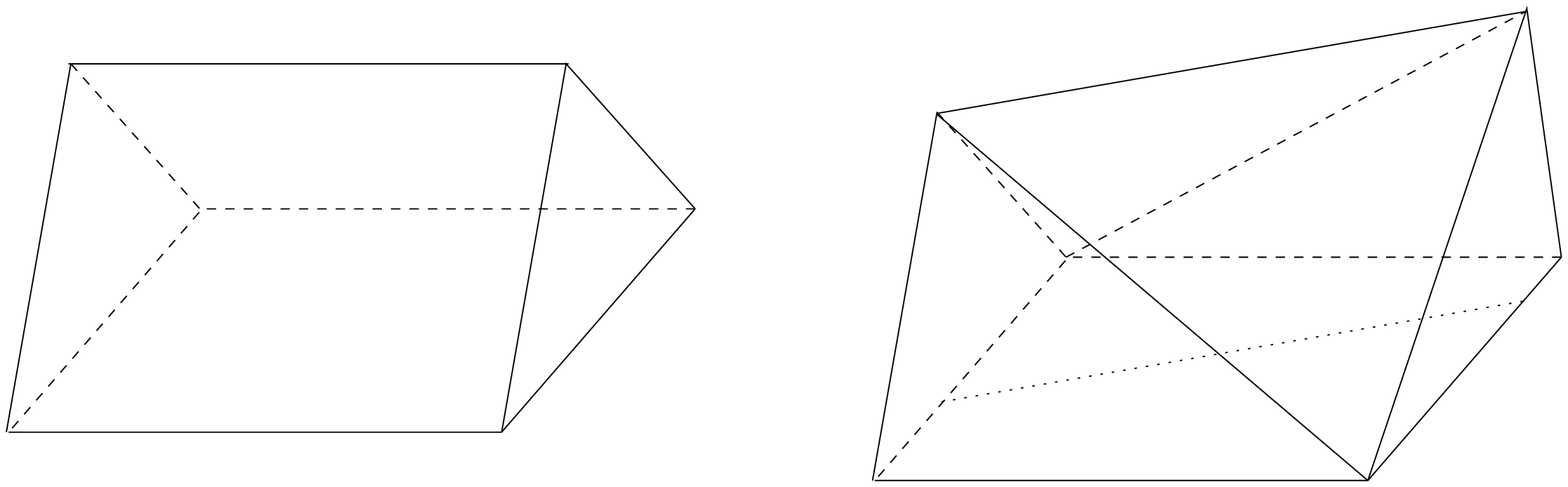}
\put(-2,0){$a$}
\put(33,0){$b$}
\put(45,16){$c$}
\put(8,16){$d$}
\put(1,27){$e$}
\put(38,27){$f$}
\put(53,-3){$a$}
\put(88,-3){$b$}
\put(101,13){$c$}
\put(63,13){$d$}
\put(58,25){$e$}
\put(99,28){$f$}
\put(58,6){$e'$}
\put(96.5,8){$f'$}
\end{overpic}
\end{center}
\caption{}\label{fig5}
\end{figure}
\end{enumerate}
\end{theorem}
\begin{proof}
We first show that if $\conv S$ has a nontriangular facet then the second case occurs.
If, on the other hand, all facets are triangular, we show that $\conv S$ must be an octahedron, and then (this being the most involved part of the proof) that some two opposite facets are parallel.

Let $P=\conv S$.
By Lemma~\ref{lemma} each nonzero point of $S-S$ is on the boundary of $P-P$.
\subsubsection*{Case I. $P$ has a nontriangular facet}
The vertex set of this facet is a planar antipodal set of more than three points, and so it must be a parallelogram $abcd$, say.
Denote the remaining two points of $S$ by $e$ and $f$.
After making an appropriate relabeling of the points and an affine transformation we have the following coordinates.
\begin{align*}
& a=(0,0,0),\; b=(1,0,0),\;  c=(1,1,0),\;  d=(0,1,0),\\
& e= (0,0,1),\;  f= (\alpha,\beta,\gamma),\; \alpha\geq\beta\geq 0,\; 0<\gamma\leq 1
\end{align*}
(We may assume $\gamma\leq 1$ after possibly interchanging $e$ and $f$.
We may assume $\alpha\geq\beta\geq 0$ after relabeling $a,b,c,d$.)

If $\beta=0$ then $e,f,a,b$ are coplanar, hence must form a parallelogram, and we obtain an affine triangular prism.
Assume then without loss of generality that $\beta>0$.
We show that this implies $\gamma=1$.

Suppose $\gamma<1$.
Consider $P-P$ and its projection onto the $xy$-plane (Figure~\ref{fig6}).
In the sequel we use the words ``above'' and ``below'' in the sense of an observer looking at $P-P$ from a point on $z$-axis with a large $z$-coordinate.
It then follows from $\gamma<1$ that $f-c$ is below the triangle with vertices $e-c, f-b, f-d$, and so $f-c\in\interior\conv\{e-c,f-b,f-d,\pm(a-c),\pm(b-d)\}$, a contradiction.
\begin{figure}
\begin{center}
\begin{overpic}[scale=0.4]{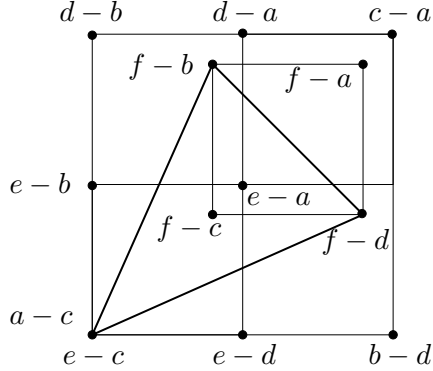}
\put(-25,5){$a-c$}
\put(-8,-8){$e-c$}
\put(40,-8){$e-d$}
\put(90,-8){$b-d$}
\put(22,33){$f-c$}
\put(75,29){$f-d$}
\put(-25,47){$e-b$}
\put(51,43){$e-a$}
\put(13,85){$f-b$}
\put(63,81){$f-a$}
\put(-9,102){$d-b$}
\put(40,102){$d-a$}
\put(90,102){$c-a$}
\end{overpic}
\end{center}
\caption{The view of $P-P$ from above}\label{fig6}
\end{figure}

Therefore, $\gamma=1$, and $f-e$ is in the $xy$-plane.
However, since the difference of any two of $a,b,c,d$ is on $\bd(P-P)$, it follows that the intersection of $P-P$ with the $xy$-plane is the square with vertices $\pm(a-c)$, $\pm(b-d)$.
Therefore, $f-e$ must be on the boundary of this square, which gives $\alpha=1$.
We now have the second type.

\subsubsection*{Case II. All facets of $P$ are triangles}
There are only two combinatorial types of $3$-polytopes with $6$ vertices and all facets triangular (by Steinitz' theorem \cite[Chapter~4]{Ziegler} it is sufficient to enumerate the $3$-connected planar triangulations on $6$ vertices).
One type (Figure~\ref{fig7}) is easily eliminated.
\begin{figure}
\begin{center}
\begin{overpic}[scale=0.2]{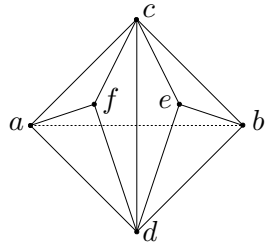}
\put(-9,48){$a$}
\put(103,48){$b$}
\put(52.5,100){$c$}
\put(52.5,-3){$d$}
\put(60,58){$e$}
\put(34,58){$f$}
\end{overpic}
\end{center}
\caption{First $3$-connected planar triangulation on $6$ vertices}\label{fig7}
\end{figure}
With the vertices labeled as shown, we apply Lemma~\ref{fivepointscor}.\eqref{one} to $S\setminus\{f\}$ to obtain that $e$ is in the half space bounded by the plane through $cd$ containing a line parallel to $ab$, opposite $a$ and $b$.
A similar argument with $S\setminus\{e\}$ gives that $f$ is also in this half space.
It follows that $\triangle cde$ and $\triangle cdf$ cannot be facets, a contradiction.

The second combinatorial type is an octahedron.
Let its diagonals be $ab$, $cd$, $ef$, say.
If each pair of diagonals is coplanar, then each such pair must be the diagonals of a parallelogram (since we then have a planar antipodal subset).
It then follows that all three diagonals are concurrent, and we obtain that $P$ is an affine regular octahedron (with any two opposite facets parallel).

In the remaining case some two diagonals are not coplanar.
It remains to show that some two opposite facets are parallel.
Without loss of generality we let $ab$ and $cd$ be noncoplanar.
After an appropriate affine transformation (mapping the vertices of the tetrahedron $abcd$ to the vertices of the cube $\{\pm1\}^3$ with an odd number of minus signs), we may assume that the $6$ points have the following coordinates (Figure~\ref{fig8}):
\begin{figure}
\begin{center}
\begin{overpic}[scale=0.25]{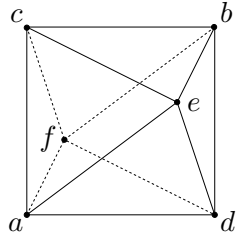}
\put(-8,-7){$a$}
\put(101,101){$b$}
\put(-7,101){$c$}
\put(101,-7){$d$}
\put(84,55){$e$}
\put(8,37){$f$}
\end{overpic}
\end{center}
\caption{Second $3$-connected planar triangulation on $6$ vertices}\label{fig8}
\end{figure}
\begin{align*}
& a=(-1,-1,-1),\; b=(1,1,-1),\;  c=(-1,1,1),\;  d=(1,-1,1),\\
& e= (\alpha,\beta,\gamma),\; \gamma>1,\;\;  f= (\alpha',\beta',\gamma'),\; \gamma'<-1,\; -\gamma'\geq\gamma.
\end{align*}
Consider the antipodal set $S\setminus\{f\}$ with convex hull $P_1$, say.
By Lemma~\ref{fivepointscor}.\eqref{two} the two planes through $a$, one parallel to $bce$ and one parallel to $bde$, both support $P_1$.
These planes have normals $(1-\beta,\alpha+\gamma,1-\beta)$ and $(\alpha+\gamma,1-\beta,1-\beta)$, respectively.
A simple calculation with inner products gives $\alpha\leq 1$ and $\beta\leq 1$.
Considering in the same way the planes through $b$ parallel to $ace$ and $ade$, we obtain $\alpha,\beta\geq -1$.
A similar argument with $P_2:=\conv(S\setminus\{e\})$ gives $-1\leq\alpha',\beta'\leq 1$.

We now consider $P-P$ and project it orthogonally onto the $xy$-plane.
The differences of pairs of $a,b,c,d$ form the $12$ vertices of a cuboctahedron that are projected onto the boundary of the square $\Sigma$ with vertices $\pm(b-a)=\pm(2,2,0)$ and $\pm(c-d)=\pm(-2,2,0)$.
Let $\Sigma_1$ ($\Sigma_2$) be the square in the $xy$-plane with vertices the projections of $e-\{a,b,c,d\}$ ($\{a,b,c,d\}-f$).
Since $-1\leq\alpha,\beta,\alpha',\beta'\leq 1$, we have $\Sigma_1,\Sigma_2\subset\Sigma$.
See Figure~\ref{fig9}.
\begin{figure}
\begin{center}
\begin{overpic}[scale=0.48]{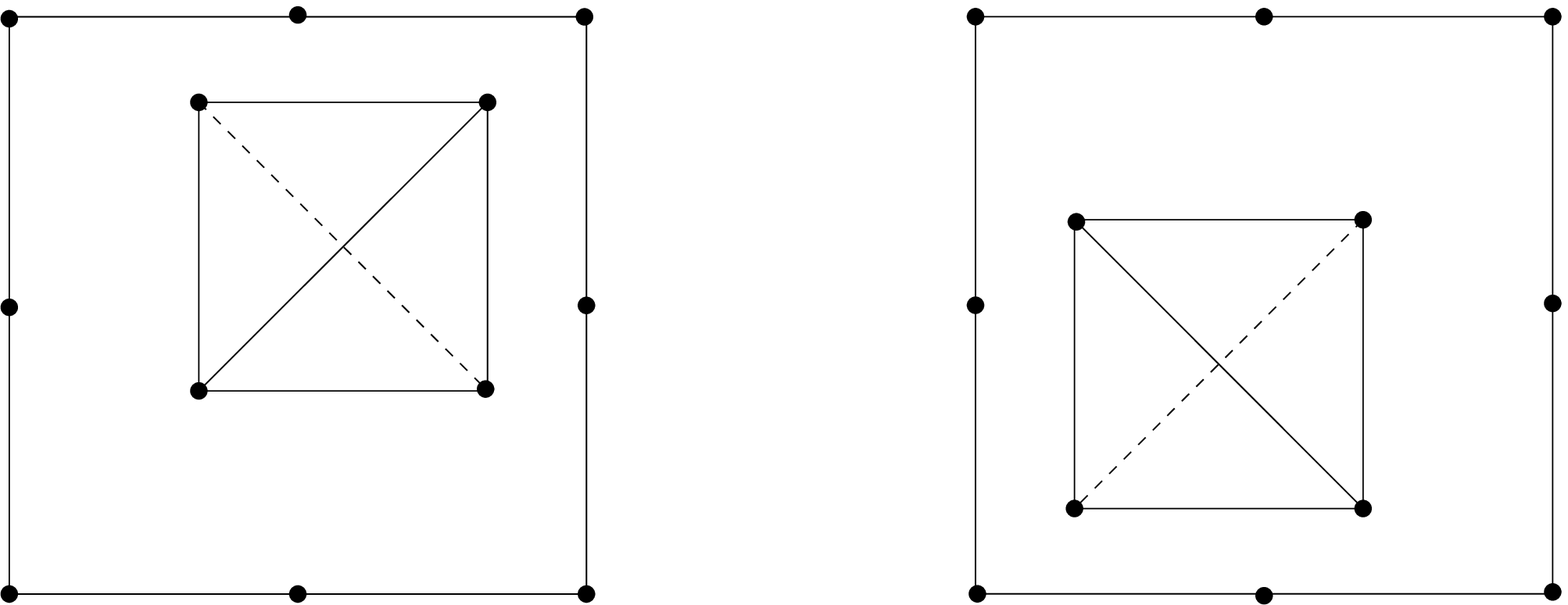}
\put(-4,-3.5){$a-b$}
\put(15,-3.5){$d-b$}
\put(33,-3.5){$d-c$}
\put(-9.5,18){$c-b$}
\put(39.5,18){$d-a$}
\put(-4,39.5){$c-d$}
\put(15,39.5){$c-a$}
\put(33,39.5){$b-a$}
\put(-5,30){$\Sigma$}
\put(7,25){$\Sigma_1$}
\put(27,34){$e-a$}
\put(8,10){$e-b$}
\put(27,10){$e-c$}
\put(8,34){$e-d$}
\put(58,-3.5){$a-b$}
\put(77,-3.5){$d-b$}
\put(95,-3.5){$d-c$}
\put(52.5,18){$c-b$}
\put(101.5,18){$d-a$}
\put(58,39.5){$c-d$}
\put(77,39.5){$c-a$}
\put(95,39.5){$b-a$}
\put(57,30){$\Sigma$}
\put(63.5,14){$\Sigma_2$}
\put(83,26){$b-f$}
\put(64,3){$a-f$}
\put(83,3){$d-f$}
\put(64,26){$c-f$}
\end{overpic}
\end{center}
\caption{$P-P$ when viewed from above}\label{fig9}
\end{figure}
In particular, $\Sigma_1$ and $\Sigma_2$ intersect, and it follows that one of the points in $e-\{a,b,c,d\}$ is projected onto $\Sigma_2$.
We now consider each of these four cases.

If $e-c$ projects onto $\Sigma_2$, then $e-c$ is below the triangle with vertices $e-a$, $e-b$, $d-f$ (Figure~\ref{fig10}).
\begin{figure}
\begin{center}
\begin{overpic}[scale=0.48]{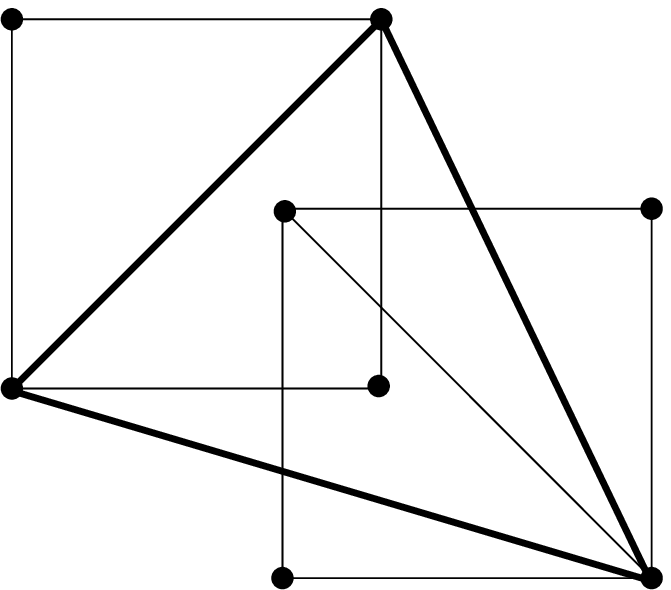}
\put(-14,75){$\Sigma_1$}
\put(47,90){$e-a$}
\put(-24,20){$e-b$}
\put(45,23){$e-c$}
\put(102,45){$\Sigma_2$}
\put(97,-8){$d-f$}
\end{overpic}
\end{center}
\caption{}\label{fig10}
\end{figure}
Since $e-c\notin\interior(P-P)$, we must have that $e-c$ projects onto the boundaries of $\Sigma_2$ and $\Sigma$, as in Figure~\ref{fig11}.
\begin{figure}
\begin{center}
\begin{overpic}[scale=0.48]{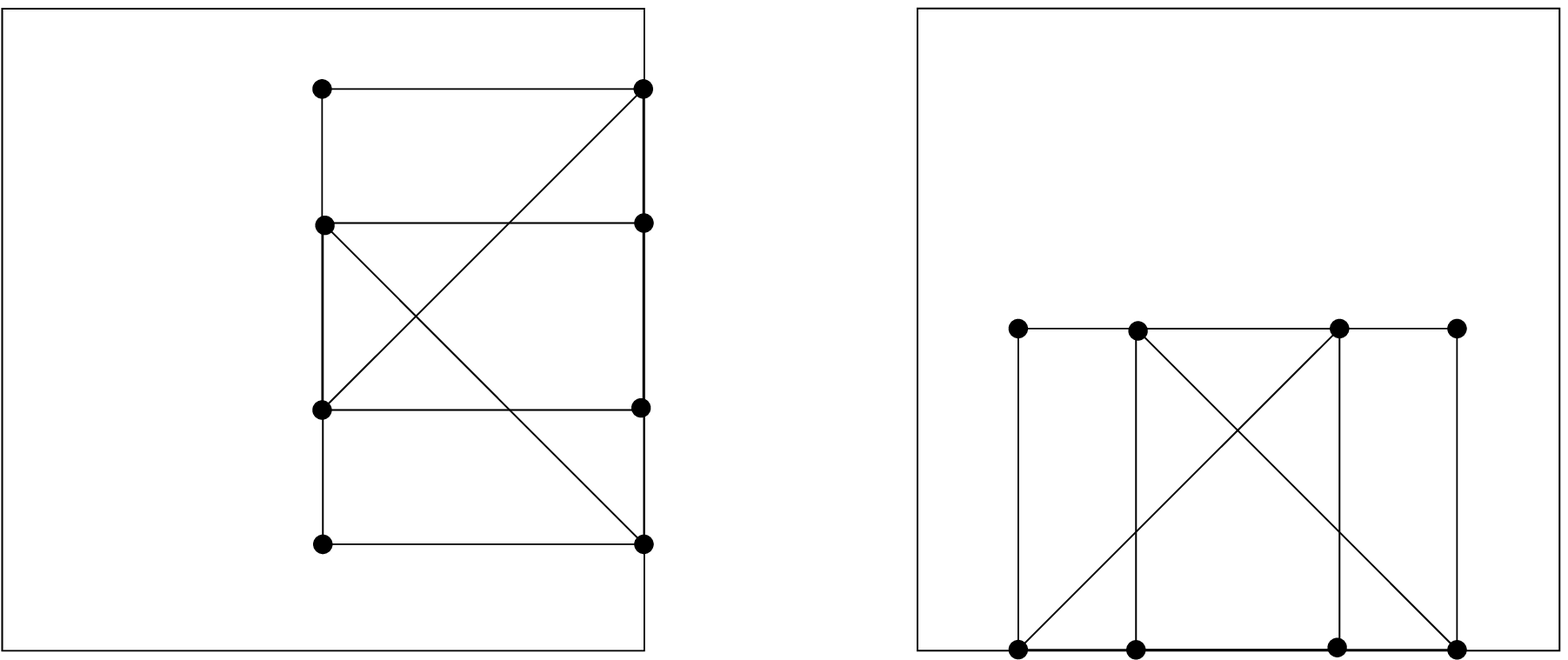}
\put(9,36){$e-d$}
\put(9,26){$c-f$}
\put(9,15){$e-b$}
\put(9,6){$a-f$}
\put(43,36){$e-a$}
\put(43,26){$b-f$}
\put(43,15){$e-c$}
\put(43,6){$d-f$}
\put(-5,30){$\Sigma$}
\put(30,38){$\Sigma_1$}
\put(30,3.5){$\Sigma_2$}
\put(68,23){$c-f$}
\put(88,23){$b-f$}
\put(60,-3){$e-b$}
\put(82,-3){$e-c$}
\put(65,30){$\Sigma$}
\put(60,13){$\Sigma_1$}
\put(94,13){$\Sigma_2$}
\end{overpic}
\end{center}
\caption{}\label{fig11}
\end{figure}
It follows that either $bde$ and $acf$ (Fig.~\ref{fig11}, left), or $ade$ and $bcf$ (Fig.~\ref{fig11}, right), are parallel, and we are finished.

A similar argument gives that if $e-d$ projects onto $\Sigma_2$, there will again be two opposite parallel facets.

Some more care is necessary with $e-a$ and $e-b$.
Suppose for instance that $e-a$ projects onto $\Sigma_2$.
If $-1-\gamma'\leq\gamma+1$, then $a-f$ is below the triangle with vertices $e-b$, $c-f$, $d-f$ (Figure~\ref{fig12}).
\begin{figure}
\begin{center}
\begin{overpic}[scale=0.58]{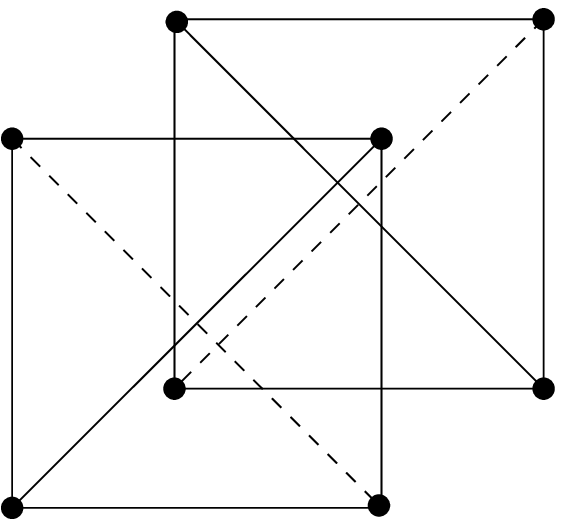}
\put(-10,72){$e-d$}
\put(52,72.5){$e-a$}
\put(-10,-8){$e-b$}
\put(53,-8){$e-c$}
\put(19,94){$c-f$}
\put(86,93){$b-f$}
\put(19,14){$a-f$}
\put(86,14){$d-f$}
\put(-13,34){$\Sigma_1$}
\put(101,58){$\Sigma_2$}
\end{overpic}
\end{center}
\caption{}\label{fig12}
\end{figure}
Since $a-f\notin\interior(P-P)$, we obtain that the projection of $a-f$ must be on the boundaries of $\Sigma_1$ and $\Sigma$, and we obtain opposite parallel facets as before.
If on the other hand $-1-\gamma'>\gamma+1$, then $e-a$ is below either $\triangle bcd-f$ or $\triangle acd-f$.
Since $e-a\notin\interior(P-P)$, we obtain that the projection of $e-a$ is on the boundaries of $\Sigma_2$ and $\Sigma$, and we again obtain parallel facets.

The case where $e-b$ projects onto $\Sigma_2$ is similar, and finishes the proof of Case II.
\end{proof}

\subsection{Seven points}
\begin{theorem}\label{sevenpoints}
Let $S$ be an antipodal set of $7$ points in $\R^3$.
Then there is a linear transformation $\fhi$ such that $\fhi(S)$ consists of the $7$ points obtained from the vertices of a cube if some two adjacent vertices of the cube are replaced by any point on the edge joining them.
\end{theorem}
\begin{proof}
The convex hull of $S$ is a $3$-polytope $P$ with $7$ vertices.
We consider various cases depending on the degrees of the vertices.

Suppose first that one of the vertices of $P$, say $a$, has degree $6$, so that it is joined by an edge to the $6$ other vertices.
Remove one of the other vertices, say $b$.
Then $S\setminus\{b\}$ will be an antipodal set of $6$ points, and in its convex hull the vertex $a$ will have a degree of $5$.
However, by Theorem~\ref{sixpoints}, no vertex can have a degree of $5$, which is a contradiction.

Suppose next that no vertex of $P$ has degree $3$.
Then all degrees are either $4$ or $5$.
Since the $1$-skeleton of $P$ is a planar graph, it has at most $15$ edges, and there are exactly two cases:
\begin{enumerate}
\item all vertices have degree $4$,
\item $5$ vertices have degree $4$, and two have degree $5$.
\end{enumerate}
There are exactly two graphs on $7$ vertices with each vertex of degree $4$, none of them planar.
In the second case there are three graphs.
In two of them the two vertices of degree $5$ are adjacent, and by removing this edge, one obtains the two graphs in which all vertices have degree $4$, which we already know to be nonplanar.
In the third graph the vertices of degree $5$ are not adjacent.
By removing one of them, the other vertex still has degree $5$, and we again obtain an antipodal set on $6$ points with a vertex of degree $5$ in its convex hull, a contradiction as before.

The only remaining possibility is that one of the vertices of $P$, say $a$, has degree $3$.
Remove a point, say $e$, that is not a neighbor of $a$.
Then the convex hull $P'$ of the antipodal set $S\setminus\{e\}$ still has $a$ as a vertex of degree $3$.
Therefore, $P'$ is not an octahedron, and by Theorem~\ref{sixpoints}, $a$ must be a vertex of a parallelogram facet $abcd$ of $P'$.
This parallelogram is also a facet of $P$.
Again using Theorem~\ref{sixpoints} we see that the remaining vertices of $P$, say $e, f, g$, are in a plane parallel to $abcd$.
Moreover, some translates of the edges $[ef], [fg], [eg]$ meet opposite sides of $abcd$.
This is only possible if one of $[ef], [fg], [eg]$ is a translate of one of the sides of $abcd$, and we obtain the conclusion.
\end{proof}

\section{Proofs of Theorems~\ref{thm2}, \ref{thm3} and \ref{thm4}}\label{lastsection}
\begin{proof}[Proof of Theorem~\ref{thm3}]
Let $S$ be an equilateral set of $6$ points at distance $1$.
Then $S$ is antipodal, and by Theorem~\ref{sixpoints} and Proposition~\ref{construction} there exist two parallel planes $\Pi_1$ and $\Pi_2$ such that $S_1=S\cap\Pi_1=\{a_1,a_2,a_3\}$ and $S_2=S\cap\Pi_2=\{b_1,b_2,b_3\}$, and the points $a_i-a_j$, $b_i-b_j$, $i\neq j$, are all in the relative boundary of their convex hull.
Also, $S_1-S_2\subset\bd B$.

Suppose that $S_1$ and $-S_2$ are translates, say with $a_i-a_j=b_j-b_i$ for all distinct $i, j$.
Then $\pm(a_i-b_j)=\pm(a_j-b_i)$ are the midpoints of the segments $\pm[a_i-b_i,a_j-b_j]$, $i<j$, and $a_i-a_j=b_i-b_j$ is the midpoint of $[a_i-b_i,b_j-a_j]$, $i\neq j$.
These $12$ segments are therefore contained in $\bd B$ and form the $1$-skeleton of an affine regular octahedron with center $o$ and vertex set $V=\{\pm(a_i-b_i):i=1,2,3\}$.
Letting $t=a_i+b_i$ (which is independent of $i$) we also have $S=\frac{1}{2}(V+t)$.

If $S_1$ and $-S_2$ are not translates, then one of the points in $S_1-S_2$ will be in the relative interior of $P:=\conv(S_1-S_2)$, which forces $P$ to be contained in $\bd B$.
Also, $P=\triangle a_1a_2a_3-\triangle b_1b_2b_3$ is a semiregular hexagon of side length $1$.

The converse is similar.
\end{proof}
\begin{proof}[Proof of Theorem~\ref{thm4}]
Let $S$ be an equilateral set of $7$ points.
By Theorem~\ref{sevenpoints}, $S$ must be as stated.
Furthermore, $\bd B$ must contain $(S-S)\setminus\{o\}$, which implies that $B$ must contain a $3$-polytope which equals some $P_\lambda$ after an appropriate linear transformation $\fhi$, and also that the planes through the facets of $[-1,1]^3$ must support $\fhi(B)$.

The converse is easy.
\end{proof}

\begin{proof}[Proof of Theorem~\ref{thm2}]
If there exists an equilateral set of size $7$, then by Theorem~\ref{thm4} the unit ball of the norm cannot be differentiable. 
\end{proof}

\section*{Acknowledgements}
We would like to thank Mohammad Ghomi for helpful explanations, the referee for suggestions leading to a better paper, and Tibor Bisztriczky for drawing our attention to the overlap between Section~\ref{antipodalsection} and \cite{BB}.

\end{document}